\documentclass[11pt]{amsart}
\usepackage{graphicx}
\usepackage{mathrsfs}
\newtheorem{thm}{Theorem}[section]
\newtheorem{cor}[thm]{Corollary}

\theoremstyle{definition}
\newtheorem{defn}[thm]{Definition}
\newtheorem{example}[thm]{Example}
\theoremstyle{remark}
\newtheorem{rem}[thm]{Remark}

\begin{document}

\title[Classes of functions and circulant matrices]{Classes of functions in the set of circulant matrices and their characterization}%
\author{Vyacheslav M. Abramov}%
\address{24 Sagan Drive, Cranbourne North, Melbourne, Victoria 3977, Australia}%

\email{vabramov126@gmail.com}%

%\thanks{}%
%\date{}%
%\dedicatory{}%
%\commby{}%
% ----------------------------------------------------------------
\begin{abstract}
In this paper, we introduce new classes of functions that extend the known classes of functions of a complex variable, such as entire functions, meromorphic functions, rational functions and polynomial functions and take values in the set of circulant matrices with complex entries. For these new classes of functions, we extend the recently obtained  characterization theorems given in [Li,   B. Q. Two elementary properties of entire functions and their applications. Amer. Math. Monthly 122(2) (2015), 169--172] and [Li, B. Q.  A characterization of rational functions. Amer. Math. Monthly. 132(3)(2025), 269--271] to an algebraic structure of circulant matrices that includes several complex variables. Our characterization theorems lead to the generalization of the variant of the fundamental theorem of algebra established recently in [Abramov, V. M. Extension of the fundamental theorem of algebra for polynomial matrix equation with circulant matrices. Amer. Math. Monthly 132(4) (2025), 356--360)].
\end{abstract}
\maketitle

\section{Introduction}  

Entire functions, meromorphic functions, rational functions, polynomial functions are classes of complex functions of complex variable that have a singificant importance in many areas of mathematics. Entire functions generalize polynomial functions, while meromorphic functions generalize rational functions. Different types of characterization theorems between the classes of functions serve an important role in undestanding the nature of these classes. Characterizations of the aforementioned classes of functions in the literature are given  in \cite[p. 143]{F}, \cite[p. 34]{GO}, \cite{H, Li1, Li2}, \cite[p. 195]{NO}.

Let $f(z)=p(z)/q(z)$ be a rational function, where $p(z)$ and $q(z)$ are polynomial functions. Without loss of generality, one can assume that $p(z)$ and $q(z)$ have no common non-constant polynomial factors. The number $\texttt{deg} \ p-\texttt{deg} \ q$, the difference between the degrees of $p(z)$ and $q(z)$, is called \textit{divisor}  of $f(z)$ (see \cite{Li2}). Li \cite{Li1, Li2} proved the following theorems.

\begin{thm}(Li \cite{Li2}.)\label{thm1}
A nonzero meromorphic function $f(z)$ in $\mathbb{C}$ is a rational function if and only if there exists the limit $\lim_{z\to\infty}zf^\prime(z)/f(z)=k$, where $k$ is the divisor of $f(z)$.
\end{thm}

\begin{thm}(Li \cite{Li1}.)\label{thm2}
Let $f(z)$ be a nonzero entire function. Then $f(z)$ has exactly $n$ zeros if and only if there is an entire function $q(z)$ such that
\[
\lim_{z\to\infty}z\left(\frac{f^\prime(z)}{f(z)}-q(z)\right)=n.
\]
\end{thm}

Both Theorems \ref{thm1} and \ref{thm2} generalize the fundamental theorem of algebra (FTA). In addition, Theorem \ref{thm2} yields the following characterization.

\begin{cor}(Li \cite{Li1}.)
A nonzero entire function $f(z)$ is a polynomial function of degree $n$ if and only if $\lim_{z\to\infty}zf^\prime(z)/f(z)=n$. 
 \end{cor} 

The analysis in \cite{Li1, Li2} strictly uses the concepts and methods of complex analysis, and hence all generalizations of the FTA given there, are related to the field $\mathbb{C}$ only. The FTA for an algebraically closed field with characteristic zero (see \cite{Al, S}) lies outside of those generalizations. 
 
The aim of the present study is to extend the aforementioned characterization theorems of \cite{Li1, Li2} of the theory of functions of a complex variable to the similar results of theory of functions of an algebraic structure that includes several complex variable. A known algebraic structure with nice properties is the set of all $d\times d$ circulant matrices (CM) with complex entries. That is, instead of complex variables $z\in\mathbb{C}$ we are going to consider $d\times d$ CM with complex entries. In fact, we are going to develop the aforementioned theorems of \cite{Li1, Li2} to the case of functions of $d$ independent complex variables in a specified way. As well, our newly obtained characterization theorems will generalize the variant of the FTA given in \cite{A}. All extentions given in the present paper are obtained in an elementary way and easily understandable. 

The rest of the paper is structured as follows. In Section \ref{s1}, we provide the necessary background information. The main results are formulated in Section \ref{s2}. The proofs are given in Section \ref{s3}. In Section \ref{s4}, we conclude the paper.

\section{Background information} \label{s1}

\subsection{} For a fixed integer $d\geq2$, the class of $d\times d$ CM with complex entries will be denoted by $CM_d$. 
A CM $\boldsymbol{X}$ has the form
\begin{equation*}
\boldsymbol{X}=\left(\begin{matrix}
x_0 &x_1 &x_2 &\cdots &x_{d-1}\\
x_{d-1} &x_0 &x_1 &\cdots &x_{d-2}\\
x_{d-2} &x_{d-1} &x_0 &\cdots &x_{d-3}\\
\vdots &\vdots &\vdots &\cdots &\vdots\\  
x_1 &x_2 &x_3 &\cdots &x_0
\end{matrix}\right).
\end{equation*}
In the sequel, it will be denoted by $\texttt{circ}(x_0, x_1, x_2,\ldots, x_{d-1})$.

The $d\times d$ circulant matrix $\boldsymbol{C}=\texttt{circ}(0, 1, 0, \ldots, 0)$ is called \textit{elementary} circulant matrix. It forms a cyclic group \cite{W1}, i.e.,  $\boldsymbol{C}^2=\texttt{circ}(0, 0, 1, \ldots, 0)$, \ldots,  $\boldsymbol{C}^{d-1}=\texttt{circ}(0, 0, 0, \ldots, 1)$, $\boldsymbol{C}^d=\texttt{circ}(1, 0, 0, \ldots, 0)=\boldsymbol{C}^0=\boldsymbol{I}$, where $\boldsymbol{I}$ is the unit matrix. We have the following canonical representation:
\begin{equation}\label{0}
\boldsymbol{X}=\sum_{i=1}^d x_{i-1}\boldsymbol{C}^{i-1}. 
\end{equation} 
The Moore-Penrose inverse (or pseudoinverse) \cite{B, M, P} of an arbitrary matrix $\boldsymbol{A}$ (invented independently by Eliakim Hastings Moore (1862--1932), Sir Roger Penrose (b. 1931) and Arne Bjerhammar (1917--2011)) is denoted by $\boldsymbol{A}^+$. If $\boldsymbol{A}$ is invertible (i.e., its determinant is not equal to zero), then $\boldsymbol{A}^+=\boldsymbol{A}^{-1}$. For a CM $\boldsymbol{X}$,  its pseudoinverse matrix $\boldsymbol{X}^+$ is circulant \cite{SB, W}.

\subsection{} Let $f(z)$ be an entire function of a complex variable $z$. It can be presented as power series
\[
f(z)=\sum_{n=0}^\infty a_nz^n
\]
for all complex $z$, i.e., it converges for all $z\in\mathbb{C}$ and converges uniformly on any compact set. 

Here we will assume that the entire function $F: CM_d\to CM_d$ inherits the properties of entire functions in $\mathbb{C}$. That is,
\[
F(\boldsymbol{Z})=\sum_{n=0}^\infty \boldsymbol{A}_n\boldsymbol{Z}^n \quad\text{for all}\ \boldsymbol{Z}\in CM_d,
\] 
where $\boldsymbol{A}_n=\texttt{circ}(a_{n,0}, a_{n,1},\ldots, a_{n,d-1})$ are matrices-coefficients, and converges uniformly on any compact set of $CM_d$. Due to canonic representation \eqref{0}, any compact set in $CM_d$ can be characterized via its counterpart in $\mathbb{C}^d$.

Meromorphic functions, rational functions, polynomial functions are also assumed to be given and taking the values in $CM_d$. For instance, a polynomial function of degree $n$ has the form $F(\boldsymbol{Z})=\boldsymbol{A}_0\boldsymbol{Z}^n+\boldsymbol{A}_1\boldsymbol{Z}^{n-1}+\ldots+\boldsymbol{A}_n$. In turn, a rational function has the presentation $F_{n,m}(\boldsymbol{Z})=P(\boldsymbol{Z})(Q(\boldsymbol{Z}))^+$, where $P(\boldsymbol{Z})$ is a polynomial function of degree $n$, and $Q(\boldsymbol{Z})$ is a polynomial function of degree $m$. The exact definition of a rational function is given later.
(The indicies $n$ and $m$ will be further omitted.) All of these definitions of functions, where all the matrices-coefficients and matrix-variable $\boldsymbol{Z}$ are CM  are correct, since CM form commutative ring \cite[p. 35]{G}.
 
Let $P(\boldsymbol{Z})=\boldsymbol{A}_0\boldsymbol{Z}^n+\boldsymbol{A}_1\boldsymbol{Z}^{n-1}+\ldots+\boldsymbol{A}_{n}$ be a polynomial function, where the matrices-coefficients and matrix-variable  $\boldsymbol{Z}$ are CM. We call the polynomial function \textit{singular}, if $\boldsymbol{A}_0$ is a singular matrix. Otherwise, it is called \textit{regular}. It follows from \cite{A} that the equation $P(\boldsymbol{Z})=\boldsymbol{O}$ with regular polynomial matrix, where $\boldsymbol{O}$ is zero matrix, satisfies the fundamental theorem of algebra. However, if $P(\boldsymbol{Z})$ is singular, a solution of $P(\boldsymbol{Z})=\boldsymbol{O}$ need not exist. For instance, if $P(\boldsymbol{Z})=\boldsymbol{E}(\boldsymbol{Z}^n+\boldsymbol{Z}^{n-1}+\ldots+\boldsymbol{Z})+\boldsymbol{I}$, where $\boldsymbol{E}$ is the matrix of units (i.e., $\boldsymbol{E}=\texttt{circ}(1, 1, \ldots, 1)$), and $\boldsymbol{I}$ is the unit matrix (i.e., $\boldsymbol{I}=\texttt{circ}(1, 0, \ldots, 0)$), it is readily seen that a solution in the form of CM of $P(\boldsymbol{Z})=\boldsymbol{O}$ does not exist. On the other hand, if $P(\boldsymbol{Z})=\boldsymbol{E}(\boldsymbol{Z}^n+\boldsymbol{Z}^{n-1}+\ldots+\boldsymbol{Z}+\boldsymbol{I})$, then the equation $P(\boldsymbol{Z})=\boldsymbol{O}$ has an infinite set of solutions in the form of CM.

Based on this, the definition of a rational function is as follows. The function $F(\boldsymbol{Z})=P(\boldsymbol{Z})(Q(\boldsymbol{Z}))^+$ is a rational function, if the equation $Q(\boldsymbol{Z})=\boldsymbol{O}$ has no more than a finite number of solutions in the form of CM. In other words, in the presentation
\begin{equation}\label{explMP}
Q(\boldsymbol{Z})=\sum_{i=0}^m\boldsymbol{B}_{i}\boldsymbol{X}^{m-i},
\end{equation}
at least one of the CM $\boldsymbol{B}_i$, $i=0,1,\ldots,m$, must be invertible. Then the total number of the solution of the equation $Q(\boldsymbol{Z})=\boldsymbol{O}$ must be not greater than $m^d$.

If the total number of solutions of the equation $Q(\boldsymbol{Z})=\boldsymbol{O}$ is not necessarily finite, the function $F(\boldsymbol{Z})$ will be called meromorphic. For instance, in the case when $Q(\boldsymbol{Z})=\boldsymbol{E}(\boldsymbol{Z}^m+\boldsymbol{Z}^{m-1}+\ldots+\boldsymbol{Z}+\boldsymbol{I})$ the function $F(\boldsymbol{Z})$ falls into the category of meromorphic but not rational, despite the degree $m\geq1$ is finite and the set of roots of $Q(\boldsymbol{Z})=\boldsymbol{O}$ is \textit{countinuously} infinite. Indeed, the equation $Q(\boldsymbol{Z})=\boldsymbol{E}(\boldsymbol{Z}^m+\boldsymbol{Z}^{m-1}+\ldots+\boldsymbol{Z}+\boldsymbol{I})=\boldsymbol{O}$ reduces to the monic polynomial equation $u^m+u^{m-1}+\ldots+1=0$, where $u=z_0+z_1+\ldots+z_{m-1}$.

\subsection{}\label{2.2+}
The reason that we consider Moore--Penrose pseudoinverse $(Q(\boldsymbol{Z}))^+$ is that the polynomial function \eqref{explMP}, when $\boldsymbol{Z}$ varies, does not need to have an inverse. For instance, suppose that $Q(\boldsymbol{Z})=\boldsymbol{Z}+\boldsymbol{E}$. Then for all $\boldsymbol{Z}$ taking the values $c\boldsymbol{E}$, $c\in\mathbb{C}$,  CM $Q(\boldsymbol{Z})$ are not invertible.

\subsection{}\label{S2.4} Let  $\omega=\mathrm{e}^{\boldsymbol{i}2\pi/d}$, where $\boldsymbol{i}=\sqrt{-1}$, and let $\omega^0$, $\omega^1$, $\omega^2$,\ldots, $\omega^{d-1}$ denote primitives of the $d$th root of unity.  Denote by $\overline{\omega}=1/\omega$ the conjugate of $\omega$. One of the most important properties of CM is that their eigenvectors are the same for all matrices from $CM_d$ (see \cite{AB, A, G, KS}). They are precisely the columns of the following matrix
\[
\boldsymbol{S}=\left(\begin{matrix} 1 &1 &1 &1 &\cdots &1\\
1 &\omega &\omega^2 &\omega^3 &\cdots &\omega^{d-1}\\
1 &\omega^2 &\omega^4 &\omega^6 &\cdots &\omega^{2d-2}\\ 
1 &\omega^3 &\omega^6 &\omega^9 &\cdots &\omega^{3d-3}\\
\vdots &\vdots &\vdots &\vdots &\cdots &\vdots\\
1 &\omega^{d-1} &\omega^{2d-2} &\omega^{3d-3} &\cdots &\omega^{(d-1)^2}
\end{matrix}\right).
\]
The matrix $\boldsymbol{S}=[s_{i,j}]$ satisfies the following properties: it is a symmetric Vandermonde matrix, and being multiplied by the factor $\sqrt{d}/d$ becomes a discrete Fourier transform matrix \cite{GG}, \cite[Chapter 2]{RY} and unitary matrix. The entries of the inverse matrix $\boldsymbol{S}^{-1}=[\widetilde{s}_{i,j}]$ satisfy the relationship $\widetilde{s}_{i,j}=d^{-1}\overline{s}_{i,j}$, where $\overline{s}_{i,j}$ denotes the conjugate of $s_{i,j}$.

The transform $\boldsymbol{S}\boldsymbol{Z}\boldsymbol{S}^{-1}=\boldsymbol{U}$ diagonalizes the CM $\boldsymbol{Z}$. The diagonal entries of the matrix $\boldsymbol{U}$ denoted by $u_i$, $i=1, 2,\ldots, d$, are defined by the formula (see \cite{A})
\begin{equation}\label{2}
u_i=\sum_{j=1}^d z_{j-1}\overline{\omega}^{(i-1)(j-1)}.
\end{equation}

\begin{example}\label{ex1}
If $\boldsymbol{Z}=\boldsymbol{E}$ and $d=p$ is a prime number, then from \eqref{2} we have $u_1=p$, and $u_i=\sum_{j=1}^p \overline{\omega}^{(i-1)j}=\sum_{j=1}^p \overline{\omega}^j=0$ for all $2\leq i\leq p$.  The equality $\sum_{j=1}^p \overline{\omega}^{(i-1)j}=\sum_{j=1}^p\overline{\omega}^j$ follows from the fact that for each $2\leq i\leq p$, and any $j_1$, $j_2$ satisfying $1\leq j_1<j_2\leq p$, the values $(i-1)j_1 \ (\texttt{mod} \ p)$ and $(i-1)j_2 \ (\texttt{mod} \ p)$ are distinct. The last is true, since otherwise $(i-1)(j_2-j_1)\ (\texttt{mod} \ p)$ must be equal to zero, which is a contradiction.
\end{example}
 
Note that if $F(\boldsymbol{Z})$ is polynomial (resp. rational, entire, meromorphic) function, then the diagonal entries of $\boldsymbol{S}F(\boldsymbol{Z})\boldsymbol{S}^{-1}$ are polynomial (resp. rational, entire, meromorphic) functions. For the details related to polynomial functions see \cite{A}.

\medskip
We will say that $\boldsymbol{Z}$ tends to infinity ($\boldsymbol{Z}\to\infty$) if  $\min_{1\leq i\leq d}|u_i|\to\infty$.

\subsection{} 

In this section, we introduce the notion of derivative of $F(\boldsymbol{Z})$.
Let $F(\boldsymbol{Z})=\texttt{circ}(f_0, f_1,\ldots, f_{d-1})$, where $f_0$, $f_1$,\ldots, $f_{d-1}$ all are the functions of $z_0$, $z_1$,\ldots, $z_{d-1}$. We define the derivative of $F(\boldsymbol{Z})$ belonging to $CM_d$ by the following way. Consider first the folowing CM of incriments
\[
\triangle [\boldsymbol{Z}]=\texttt{circ}(\triangle z_0, \triangle z_1,\ldots, \triangle z_{d-1}),
\]
assuming that $\triangle [\boldsymbol{Z}]$ is invertible. Then, it is convenient to assume that $\triangle z_i=\delta c_i$, where $\texttt{circ}(c_0, c_1,\ldots, c_{d-1})$ is an invertible CM, and $\delta$ is a variable parameter.

The limit
\[
\begin{aligned}
&\lim_{\delta\to0}\frac{ f_i(z_0, z_1,\ldots, z_j+\delta c_j,\ldots, z_{d-1})-f_i(z_0,z_1,\ldots,z_{d-1})}{\delta c_j}\\  %&=\delta_{i,j}(z_0,z_1,\ldots,z_{d-1}), \quad i,j=0,1,\ldots, d-1,
\end{aligned}
\]
if it exists, denotes the partial derivative $\partial f_i/\partial z_j$.

Following this, introduce the formal notation
\begin{eqnarray*}
\mathrm{d} [F(\boldsymbol{Z})]&=&\texttt{circ}(\partial f_0, \partial f_1,\ldots, \partial f_{d-1}),\\
\mathrm{d}[\boldsymbol{Z}]&=& \texttt{circ}(\partial z_0, \partial z_1,\ldots, \partial z_{d-1}).
\end{eqnarray*}
Then, by the derivative of $F(\boldsymbol{Z})$ we mean
\begin{equation}\label{12}
\mathscr{D}[F(\boldsymbol{Z})]=\mathrm{d} [F(\boldsymbol{Z})](\mathrm{d}[\boldsymbol{Z}])^{-1},
\end{equation}
assuming that the partial derivatives $\partial f_i/\partial z_j$ for all $i,j=0, 1,\ldots, d-1$ exist.

The definition of the derivative seems to be too formal. However, if $(\mathrm{d}[\boldsymbol{Z}])^{-1}$ is presented in the form
\[
(\mathrm{d}[\boldsymbol{Z}])^{-1}=\texttt{circ}\left(\frac{a_0}{\partial z_0}, \frac{a_1}{\partial z_1},\ldots, \frac{a_{d-1}}{\partial z_{d-1}}\right),
\] 
for some constants $a_0$, $a_1$,\ldots, $a_{d-1}$ all depending of $\partial z_0$, $\partial z_1$,\ldots, $\partial z_{d-1}$ (that is generally true), then the entries of $\mathscr{D}[F(\boldsymbol{Z})]$ can be written in an explicit way. For instance, for the first entry of the matrix we have
\[
a_0\frac{\partial f_0}{\partial z_0}+a_{d-1}\frac{\partial f_1}{\partial z_{d-1}}+\ldots+a_1\frac{\partial f_{d-1}}{\partial z_1}.
\]

The conjugating of $\mathscr{D}[F(\boldsymbol{Z})]$ by $\boldsymbol{S}$ yields the diagonal matrix
\[
\begin{aligned}
\boldsymbol{S}\mathscr{D}[F(\boldsymbol{Z})]\boldsymbol{S}^{-1}&=(\boldsymbol{S}\mathrm{d}[F(\boldsymbol{Z})]\boldsymbol{S}^{-1})(\boldsymbol{S}(\mathrm{d}[\boldsymbol{Z}])^{-1}\boldsymbol{S}^{-1})\\
&=(\boldsymbol{S}\mathrm{d}[F(\boldsymbol{Z})]\boldsymbol{S}^{-1})(\boldsymbol{S}\mathrm{d}[\boldsymbol{Z}]\boldsymbol{S}^{-1})^{-1},
\end{aligned}
\]
the $i$th diagonal entry of which according to \eqref{2} is equal to
\begin{equation}\label{10}
\frac{\sum_{j=1}^d \partial f_{j-1}\overline{\omega}^{(i-1)(j-1)}}{\sum_{j=1}^d \partial z_{j-1}\overline{\omega}^{(i-1)(j-1)}}
=\frac{\mathrm{d}\left(\sum_{j=1}^d  f_{j-1}\overline{\omega}^{(i-1)(j-1)}\right)}{\mathrm{d}\left(\sum_{j=1}^d  z_{j-1}\overline{\omega}^{(i-1)(j-1)}\right)},
\end{equation}
where the expression on the left-hand side of \eqref{10} characterizes the sum of the directional derivatives that after the change of the variables becomes the derivative of the function of one variable. The denominator in \eqref{10} is nonzero due to the definition of $\mathrm{d}[\boldsymbol{Z}]$ that states that CM $\mathrm{d}[\boldsymbol{Z}]$ is invertible.

Notice that the expression $\sum_{j=1}^d  f_{j-1}\overline{\omega}^{(i-1)(j-1)}$ in the numerator of the right-hand side of \eqref{10} is the $i$th diagonal entry of $\boldsymbol{S}F(\boldsymbol{Z})\boldsymbol{S}^{-1}$. 
 
\begin{example}\label{ex2}
Let $F(\boldsymbol{Z})=\boldsymbol{Z}^n$, where $n$ is a positive integer number. From the right-hand side of \eqref{10} for the $i$th diagonal entry of $\boldsymbol{S}\mathscr{D}[\boldsymbol{Z}^n]\boldsymbol{S}^{-1}$, $i=1,2,\ldots,d$, we have
\[
\frac{\mathrm{d}\left(\sum_{j=1}^d  z_{j-1}\overline{\omega}^{(i-1)(j-1)}\right)^n}{\mathrm{d}\left(\sum_{j=1}^d  z_{j-1}\overline{\omega}^{(i-1)(j-1)}\right)}=n\left(\sum_{j=1}^d  z_{j-1}\overline{\omega}^{(i-1)(j-1)}\right)^{n-1},
\]
that coincides with the corresponding $i$th diagonal entry of $\boldsymbol{S}(n\boldsymbol{Z}^{n-1})\boldsymbol{S}^{-1}$.
\end{example}

\begin{example}\label{ex3}
Let $F(\boldsymbol{Z})=P(\boldsymbol{Z})(Q(\boldsymbol{Z}))^+$, where $P(\boldsymbol{Z})$ and $Q(\boldsymbol{Z})$ are polynomial functions of degrees $n$ and $m$, respectively. Denote 
\[
u_i=\sum_{j=1}^d  z_{j-1}\overline{\omega}^{(i-1)(j-1)}.
\]

Then, for the $i$th diagonal entry of $\boldsymbol{S}\mathscr{D}[P(\boldsymbol{Z})(Q(\boldsymbol{Z}))^+]\boldsymbol{S}^{-1}$, $i=1,2,\ldots,d$, we have
\[
\frac{\mathrm{d}(P_i(u_i)(Q_i(u_i))^{-1})}{\mathrm{d}u_i}=\frac{P_i^\prime(u_i)Q_i(u_i)-Q_i^\prime(u_i)P_i(u_i)}{(Q_i(u_i))^2}.
\]
where $P_{i}(u_i)$ and $Q_{i}(u_i)$ denote the polynomial functions that appear in the $i$th diagonal entry of  $\boldsymbol{S}P(\boldsymbol{Z})\boldsymbol{S}^{-1}$ and $\boldsymbol{S}Q(\boldsymbol{Z})\boldsymbol{S}^{-1}$, respectively.
\end{example}

\section{Main results}\label{s2}

We start from the following definition. 
 
\begin{defn}\label{dfn1}
Let $F(\boldsymbol{Z})=P(\boldsymbol{Z})(Q(\boldsymbol{Z}))^+$ be a rational function in $CM_d$, where $P(\boldsymbol{Z})$ and $Q(\boldsymbol{Z})$ are corresponding regular polynomial functions of degrees $n$ and $m$. 
The value $k=n-m$ is called \textit{divisor}.
\end{defn}

\begin{rem} 
Definition \ref{dfn1} is equivalent to that given in \cite{Li2}.
\end{rem}

\begin{thm}\label{thm3}
A nonzero meromorphic function $F(\boldsymbol{Z})\in CM_d$ is a rational function if and only if
\begin{equation}\label{6}
\lim_{\boldsymbol{Z}\to\infty}\boldsymbol{S}\left(\boldsymbol{Z}\mathscr{D}[F(\boldsymbol{Z})](F(\boldsymbol{Z}))^{+}\right)\boldsymbol{S}^{-1}=\left(\begin{matrix}k_1 &0 &\cdots &0\\
0 &k_2 &\cdots &0\\
\vdots &\vdots &\ddots &\vdots\\
0 &0 &\cdots &k_d\end{matrix}\right),
\end{equation}
where $-m\leq k_i\leq n$, $i=1,2,\ldots, d$.

In addition, if the polynomials $P(\boldsymbol{Z})$ and $Q(\boldsymbol{Z})$ in the presentation of the rational function both are regular, then the right-hand side of \eqref{6} is equal to $k\boldsymbol{I}$.
\end{thm}

\begin{thm}\label{thm4}
Let $F(\boldsymbol{Z})\in CM_d$ be a nonzero entire function. Then the number of solutions of $F(\boldsymbol{Z})=\boldsymbol{O}$  does not exceed $n^d$ if and only if there is an entire function $Q(\boldsymbol{Z})\in CM_d$ such that
\[
\lim_{\boldsymbol{Z}\to\infty}\boldsymbol{S}\left(\boldsymbol{Z}\mathscr{D}[F(\boldsymbol{Z})](F(\boldsymbol{Z}))^{+}-\boldsymbol{Z}Q(\boldsymbol{Z})\right)\boldsymbol{S}^{-1}=n\boldsymbol{I}.
\]
\end{thm}

\begin{cor}\label{cor2}
A nonzero entire function $F(\boldsymbol{Z})\in CM_d$ is a regular polynomial function of degree $n$ if and only if $\lim_{\boldsymbol{Z}\to\infty}\boldsymbol{S}(\boldsymbol{Z}\mathscr{D}[F(\boldsymbol{Z})](F(\boldsymbol{Z}))^{+}\boldsymbol{S}^{-1}=n\boldsymbol{I}$. 
 \end{cor}  

\section{Proofs}\label{s3}
\subsection{Proof of Theorem \ref{thm3}.}
We have: 
\begin{equation}\label{3}
\begin{aligned}
&\underbrace{\boldsymbol{S}\left(\boldsymbol{Z}\mathscr{D}[F(\boldsymbol{Z})](F(\boldsymbol{Z}))^{+}\right)\boldsymbol{S}^{-1}}_{\text{diagonal matrix}}\\
&=\underbrace{\left(\boldsymbol{S}\boldsymbol{Z}\boldsymbol{S}^{-1}\right)}_{\text{diagonal matrix 1}}\underbrace{\left(\boldsymbol{S}\mathscr{D}[F(\boldsymbol{Z})]\boldsymbol{S}^{-1}\right)}_{\text{diagonal matrix 2}}\underbrace{\left(\boldsymbol{S}(F(\boldsymbol{Z}))^{+})\boldsymbol{S}^{-1}\right)}_{\text{diagonal matrix 3}}.
\end{aligned}
\end{equation}

The entries of diagonal matrix 1 of the right-hand side of \eqref{3} are given by \eqref{2}. The entries of diagonal matrix 2 are given by \eqref{10}. Denote the denominator of \eqref{10} by $\mathrm{d}u_i$. Then, with appropriate change of the variables, the expression in \eqref{10} for the given case can be rewritten as $\mathrm{d}F_i(u_i)/\mathrm{d}u_i:=F_i^\prime(u_i)$  for all $u_i$ where the derivative exists. Such possibility of the change of variable for the derivative is supported by Examples \ref{ex2} and \ref{ex3}.

The possibility of the aforementioned change of variables for the function itself for diagonal matrix 3 is confirmed in an elementary way as well. For instance, if $F(\boldsymbol{Z})$ is a rational function, i.e. $F(\boldsymbol{Z})=P(\boldsymbol{Z})(Q(\boldsymbol{Z}))^+$, where $P(\boldsymbol{Z})$ and $Q(\boldsymbol{Z})$ are polynomial functions, $P(\boldsymbol{Z})=\boldsymbol{O}$ and   $Q(\boldsymbol{Z})=\boldsymbol{O}$ both having a finite number of solutions, we obtain
\[
\boldsymbol{S}(F(\boldsymbol{Z}))^+\boldsymbol{S}^{-1}=(\boldsymbol{S}Q(\boldsymbol{Z})\boldsymbol{S}^{-1})(\boldsymbol{S}(P(\boldsymbol{Z}))^+\boldsymbol{S}^{-1}).
\]

But 
\[
\boldsymbol{S}Q(\boldsymbol{Z})\boldsymbol{S}^{-1}=\boldsymbol{S}\left(\sum_{i=0}^m\boldsymbol{B}_i\boldsymbol{Z}^{m-i}\right)\boldsymbol{S}^{-1}=\sum_{i=0}^n(\boldsymbol{S}\boldsymbol{B}_i\boldsymbol{S}^{-1})(\boldsymbol{S}\boldsymbol{Z}\boldsymbol{S}^{-1})^{m-i},
\]
where $\boldsymbol{B}_i$, $i=0, 1,\ldots, m$, are the matrices-coefficients of the polynomial function $Q(\boldsymbol{Z})$.

As well, taking into account that $\boldsymbol{S}(P(\boldsymbol{Z}))^+\boldsymbol{S}^{-1}=(\boldsymbol{S}P(\boldsymbol{Z})\boldsymbol{S}^{-1})^+$, in the similar way we have
\[
\begin{aligned}
(\boldsymbol{S}P(\boldsymbol{Z})\boldsymbol{S}^{-1})^+&=\left(\boldsymbol{S}\left(\sum_{i=0}^n\boldsymbol{A}_i\boldsymbol{Z}^{n-i}\right)\boldsymbol{S}^{-1}\right)^+\\ &=\left(\sum_{i=0}^n(\boldsymbol{S}\boldsymbol{A}_i\boldsymbol{S}^{-1})(\boldsymbol{S}\boldsymbol{Z}\boldsymbol{S}^{-1})^{n-i}\right)^+,
\end{aligned}
\]
where $\boldsymbol{A}_i$, $i=0, 1,\ldots, n$, are the matrices-coefficients of the polynomial function $P(\boldsymbol{Z})$.

Note that the entries of 
\begin{equation}\label{P+}
\left(\sum_{i=0}^n(\boldsymbol{S}\boldsymbol{A}_i\boldsymbol{S}^{-1})(\boldsymbol{S}\boldsymbol{Z}\boldsymbol{S}^{-1})^{n-i}\right)^+
\end{equation} 
that appear in the denominator may take zero value, and the fractions in the diagonal entries can be undefined. However, if $F(\boldsymbol{Z})$ is a rational function, then for large $\boldsymbol{Z}$ all diagonal entries of \eqref{P+} are nonzero (due to the definition of large $\boldsymbol{Z}$ given at the end of Section \ref{S2.4}).

The term $\boldsymbol{S}\boldsymbol{Z}\boldsymbol{S}^{-1}$ is fully separated, that enables us to change the required variables that appear in the entries of diagonal matrix 3. (The explicit presentation of $\boldsymbol{S}\boldsymbol{Z}\boldsymbol{S}^{-1}$ is given by \eqref{2}.)

So, the $i$th diagonal entry on the left-hand side of \eqref{3} is presented in the form
\begin{equation}\label{ithdiagentry}
\frac{u_iF_i^\prime(u_i)}{F_i(u_i)},
\end{equation}
and is defined for all large $u_i$, if $F_i(u_i)$ is a rational function.

If $F(\boldsymbol{Z})$ is a meromorphic function and not a rational function, the arguments are similar to the above,   the  $i$th diagonal entry given by \eqref{ithdiagentry}  as well, but the function $F_i(u_i)$ is meromorphic and not rational.

Now, according to Theorem \ref{thm1} any $i$th diagonal entry on the left-hand side of \eqref{3} must converge to $k_i$ if and only if $
F_i(u_i)$ is rational of divisor $k_i$. But if $F(\boldsymbol{Z})$ is rational, then all of $F_i(u_i)$ must be rational, and for large values of $u_i$ we have $F_i(u_i)=P_i(u_i)/Q_i(u_i)$. If both polynomials $P(\boldsymbol{Z})$ and $Q(\boldsymbol{Z})$ are regular, then the degrees of $P_i(u_i)$ and $Q_i(u_i)$ must be $n$ and $m$, respectively, and hence divisor of $F_i(u_i)$ is equal to $k=n-m$. Otherwise, if at least one of the polynomials $P(\boldsymbol{Z})$ or $Q(\boldsymbol{Z})$ is singular, then the degrees of $P_i(u_i)$ or $Q_i(u_i)$ may be less than $n$ or, respectively, less than $m$. Indeed, following Example \ref{ex1} for $\boldsymbol{E}\boldsymbol{Z}^n$ with the $p\times p$ matrices, where $p$ is a prime number, we have as follows:
\[
\boldsymbol{S}(\boldsymbol{E}\boldsymbol{Z}^n)\boldsymbol{S}^{-1}=(\boldsymbol{S}\boldsymbol{E}\boldsymbol{S}^{-1})(\boldsymbol{S}\boldsymbol{Z}\boldsymbol{S}^{-1})^n:=\boldsymbol{D}.
\]
The resulting diagonal matrix $\boldsymbol{D}$ satisfies the properties: the first diagonal element is equal to $d(z_0+z_1+\ldots+z_{p-1})^n$, and all the other diagonal elements are zeros.

Hence, the divisor of $F_i(u_i)$ denoted by $k_i$ must satisfy the inequality $-m\leq k_i\leq n$. The proof is completed.

\begin{rem}\label{rem3}
If $F(\boldsymbol{Z})=P(\boldsymbol{Z})$ is a regular polynomial function, then the limit in \eqref{6} must be equal to $n\boldsymbol{I}$. It follows from the proof given in \cite{Li2} that the fact that $\lim_{u_i\to\infty}u_iF_i^\prime(u_i)/F_i(u_i)=n$ means that the monic polynomial equation $F_i(u_i)=0$ has exactly $n$ roots counting multiplicity. The last means that the polynomial equation $F(\boldsymbol{Z})=\boldsymbol{O}$ has at least one and at most $n^d$ solutions being CM. So, Theorem \ref{thm3} extends \cite[Theorem 2]{A}.
\end{rem}

\subsection{Proof of Theorem \ref{thm4}.} We have:
\begin{equation}\label{5}
\begin{aligned}
&\underbrace{\boldsymbol{S}\left(\boldsymbol{Z}\mathscr{D}[F(\boldsymbol{Z})](F(\boldsymbol{Z}))^{+}-\boldsymbol{Z}Q(\boldsymbol{Z})\right)\boldsymbol{S}^{-1}}_{\text{diagonal matrix}}\\
&=\underbrace{\boldsymbol{S}\left(\boldsymbol{Z}\mathscr{D}[F(\boldsymbol{Z})](F(\boldsymbol{Z}))^{+})\right)\boldsymbol{S}^{-1}}_{\text{diagonal matrix 1}}-\underbrace{\boldsymbol{S}\left(\boldsymbol{Z}Q(\boldsymbol{Z})\right)\boldsymbol{S}^{-1}}_{\text{diagonal matrix 2}}.
\end{aligned}
\end{equation} 
The arguments for diagonal matrix 1 are the same as in the proof of Theorem \ref{thm3}, and the arguments for diagonal matrix 2 are standard, since
\[
\boldsymbol{S}\left(\boldsymbol{Z}Q(\boldsymbol{Z})\right)\boldsymbol{S}^{-1}=(\boldsymbol{S}\boldsymbol{Z}\boldsymbol{S}^{-1})\left(\boldsymbol{S}Q(\boldsymbol{Z})\boldsymbol{S}^{-1}\right).
\]
Therefore we can apply Theorem \ref{thm2}. According to this theorem any $i$th diagonal entry of the matrix on the left-hand side of \eqref{5} must be equal to $n$ if and only if an appropriate entire function can be chozen. The proof is completed.

\begin{rem}
Theorem \ref{thm4} extends \cite[Theorem 2]{A} as well. Indeed, if $F(\boldsymbol{Z})$ is a regular polynomial function of degree $n$, then $Q(\boldsymbol{Z})$ can be taken to equal to zero matrix, and we arrive at  the considered case in Remark \ref{rem3}. 
\end{rem}

\medskip

The proof of Corollary \ref{cor2} is similar.
 
\section{Concluding remarks}\label{s4} 
In the present paper, we extended certain characterization theorems of the theory of functions of a complex variable to the similar theory of function of an algebraic structure of CM. The methods used in the present paper can be further developed for new problems, say, in the theory of differential equations.

\subsection*{Data availability statement}
Data sharing not applicable to this article as no datasets were generated or analysed during the current study.

\subsection*{Disclosures and declarations}
No conflict of interest was reported by the author.

\subsection*{Authorship clarified}
The author, who conducted this work is not affiliated to any institution and solely contributed to this work. The author has full responsibility on the content of the paper.

\end{document}